\documentclass{article}

\usepackage{delarray,verbatim,enumerate,a4wide}
\usepackage{amsmath,amsthm,amstext,amsbsy,amssymb,amsfonts,amscd}
\usepackage[all]{xy}
\usepackage[frenchb]{babel}
\usepackage[T1]{fontenc}
\usepackage[utf8]{inputenc}

\usepackage{scalerel}

\usepackage[small]{titlesec}

\usepackage{color}
\definecolor{vert}{rgb}{0.1,0.4,0.2}
\usepackage[colorlinks=true,linkcolor=blue,citecolor=vert]{hyperref}

\usepackage{calligra}
\DeclareFontShape{T1}{calligra}{m}{n}{<->s*[0.95]callig15}{}
\DeclareMathAlphabet{\mathscr}{T1}{calligra}{m}{n}

\newtheorem{Th}{Théorème}[]

\newtheorem{Prop}[Th]{Proposition}
\newtheorem{Cor}[Th]{Corollaire}
\newtheorem{Conj}[Th]{Conjecture}
\newtheorem{Sco}[Th]{Scolie}

\newtheorem*{Th*}{Théorème}
\newtheorem*{LemA}{Lemme A}
\newtheorem*{LemB}{Lemme B}
\newtheorem*{LemC}{Lemme C}
\newtheorem*{LemD}{Lemme D}

\def\Preuve{\noindent {\it Preuve.~}}

	\def\ZZ{\mathbb Z}		
\def\F2{\mathbb{F}_2}	\def\Z2{\mathbb{Z}_2}		
\def\Zl{\mathbb{Z}_\ell} 			

 				\def\U{\mathcal  U}	\def\F{\mathcal  F}
\def\J{\mathcal  J}  	\def\C{\mathcal  C}	\def\R{\mathcal  R}	
 	  	\def\Cl{\mathcal  C\!\ell}	\def\V{\mathcal  V}

\def\G{\mathscr G\,}

		\def\p{{\mathfrak p}}				\def\a{{\mathfrak a}}		
						 \def\f{{\mathfrak f}}

\def\wi{\widetilde}			\def\End{\operatorname{End}}
		\def\Tr{\operatorname{Tr}}	\def\res{\operatorname{res}}
	\def\deg{\operatorname{deg}}	\def\Ver{\operatorname{Ver}}
\def\Gal{\operatorname{Gal}}			\def\Im{\operatorname{Im}}
\def\Ker{\operatorname{Ker}}	\def\Coker{\operatorname{Coker}}	\def\Hom{\operatorname{Hom}}

\newcommand\scale[2]{\vstretch{#1}{\hstretch{#1}{#2}}}

\newcommand\si[1]{\scale{.8}{#1}}

\makeatletter
\newcommand*\wt[2][0.2ex]{%
        \begingroup
        \mathchoice{\wt@helper{#1}{#2}{\displaystyle}{\textfont}}
                   {\wt@helper{#1}{#2}{\textstyle}{\textfont}}
                   {\wt@helper{#1}{#2}{\scriptstyle}{\scriptfont}}
                   {\wt@helper{#1}{#2}{\scriptscriptstyle}{\scriptscriptfont}}%
        \endgroup
        #2%
}
\newcommand*\wt@helper[4]{%
        \def\currentfont{\the#41}%
        \def\currentskewchar{\char\the\skewchar\currentfont}%
        \setbox\tw@\hbox{\currentfont$#2$\currentskewchar}%
        \dimen@ii\wd\tw@
        \setbox\tw@\hbox{\currentfont$#2${}\currentskewchar}%
        \advance\dimen@ii-\wd\tw@
        \rlap{\raisebox{-#1}{$\m@th#3\kern\dimen@ii\widetilde{\phantom{#2}}$}}%
}
\makeatother

\begin{document}

\title{\Large\bf L'état actuel du problème de la capitulation\footnote{Sém. Théor. Nombres Bordeaux 1987--1988, exp. n$^\circ 17$}}

\author{ Jean-François {\sc Jaulent} }
\date{}
\maketitle
\bigskip\bigskip

{\footnotesize
\noindent{\bf Résumé.} Le texte qui suit est une présentation synthétique de l'état des connaissances sur le problème de la capitulation pour les groupes de classes des corps de nombres, peu avant la démonstration par Suzuki \cite{Suz} de la conjecture principale sur cette question.}\smallskip

{\footnotesize
\noindent{\bf Abstract.} The text which follows is a synthetic presentation of the state of the knowledge about the capitulation for the class-groups of numbers fields, shortly before the demonstration by Suzuki \cite{Suz} of the main conjecture on this question.}
\smallskip

\tableofcontents

\section*{Introduction}
\addcontentsline{toc}{section}{Introduction}

On dit depuis Hilbert qu'un idéal de l'anneau des entiers d'un corps de nombres $k$ capitule dans une extension finie $K$ de $k$, lorsqu'il devient principal par extension des scalaires à l'anneau des entiers de $K$. Le problème de la capitulation consiste donc précisément à décrire la partie du groupe des classes d'idéaux de $k$ qui est représentée par les classes des idéaux de $k$ qui capitulent dans $K$, lorsque $K$ est une extension abélienne non ramifiée de $k$. Comme nous le verrons plus loin, cette hypothèse de non-ramification n'est pas vraiment une restriction: pour les extensions abéliennes de conducteur donné, le problème de la capitulation se pose de façon naturelle en termes de classes de rayons. On peut aussi le poser bien entendu en termes de classes d'idèles, ce qui se révèle techniquement commode à l'occasion, mais peu éclairant en fin de compte pour la question qui nous préoccupe.\smallskip

Quelle que soit sa formulation cependant, y compris la plus élémentaire, la question de la capitulation reste aujourd'hui encore l'un des aspects les moins connus de l'arithmétique des extensions abéliennes. Cela tient, je crois, au fait que la Théorie du corps de classes, qui est l'information la plus forte que  nous possédions sur l'arithmétique de ces extensions, s'est élaborée paradoxalement non pas, comme on a pu l'espérer avant les travaux de Takagi, à partir des propriétés de l'homomorphisme d'extension, qui sont au cœur du problème de la capitulation, mais, tout au contraire, autour de celles, duales, de l'application norme. Arrêtons nous un instant sur ce point:
étant donnée une extension quelconque $K/k$ de corps de nombres, il existe deux applications naturelles entre les groupes de classes d' idéaux de $K$ et de $k$.
\begin{displaymath}
\xymatrix{\;Cl_K \ar@<1ex>[d]^{N_{K/k}}\\ Cl_k  \!\!\ar@{->}[u]^{j_{K/k}} }
\end{displaymath}
\begin{itemize}
\item La première, ascendante, est {\em l'homomorphisme d'extension} $j_{K/k}$ qui est induit par l'extension des idéaux de l'anneau des entiers de $k$ à celui de $K$.
\item La seconde, descendante, est la {\em norme arithmétique} $N_{K/k}$.
\end{itemize}
Bien entendu, ces deux applications ne sont pas indépendantes puisqu'on a trivialement:\smallskip

\centerline{$N_{K/k}\circ j_{K/k}=[K/k] $ \qquad et \qquad $j_{K/k}\circ N_{K/k}=\nu_{K/k}$,}\smallskip

\noindent i.e. l'exponentiation par le degré de l'extension, à gauche; et la {\em norme algébrique}, à droite; c'est-à-dire l'exponentiation symbolique par l'élément $\nu_{K/k}=\sum_{g\in G}g$ de l'algèbre $\ZZ[G]$, lorsque l'extension $K/k$ est galoisienne de groupe $G$.
\smallskip

Maintenant, le sous-groupe $Cap_{K/k}$ des classes de $Cl_k$ qui capitulent dans $Cl_K$ n'est rien d'autre que le noyau du morphisme $j_{K/k}$. Mais c'est du conoyau de $N_{K/k}$ et non du noyau de $j_{K/k}$ dont nous parle la théorie du corps de classes. Cela ne veut pas dire que cette théorie ne nous apporte aucune information sur la capitulation, car nombre d'avancées significatives, comme le théorème d'Artin-Furtwängler ou celui de Tannaka-Terada ont été obtenues via le corps de classes. Simplement, ces informations sont essentiellement indirectes, ce qui explique que bien des aspects du phénomène de la capitulation nous demeurent mystérieux.
\medskip

Cela dit, le plan de ce rapport est le suivant.\smallskip

Dans une première partie, nous présentons ce que nous appelons l'interprétation arithmétique du théorème de l'idéal principal, qui consiste à relier l'étude de la capitulation à celle de la cohomologie des unités : c'est l'approche la plus ancienne, celle inaugurée par Hilbert dans son célèbre théorème 94.
\smallskip

 Dans une seconde partie, nous nous attachons à décrire l'approche algébrique du problème qui permet, via les isomorphismes du corps de classes de ramener le théorème de l'idéal principal dans certaines situations bien précises à un problème ment algébrique de théorie des groupes : c'est la méthode suivie pour les théorèmes d'Artin-Furwängler et de Tannaka-Terada évoqués plus haut.

Enfin , dans une dernière partie, nous présentons quelques résultats récents sur le problème de la capitulation ; en particulier ceux obtenus par Miyake.
\smallskip

Dans chacune des trois sections, nous avons essayé de donner des preuves aussi concises et directes que possible des résultats présentés, notamment chaque fois qu'une démonstration immédiatement accessible au lecteur non spécialiste nous sait manquer. Si donc les résultats énoncés ici sont, pour l'essentiel, bien connus des spécialistes, certaines approches peuvent être regardées, elles, comme originales.

\section{L'approche arithmétique}

\subsection{Le théorème 94 de Hilbert et son interprétation cohomologique}

Le résultat le plus ancien sur la capitulation est celui donné par D. Hilbert dans son traité sur les corps de nombres algébriques \cite{Hil}, à propos des extensions cycliques relatives de degré premier impair, et qu'il énonce comme suit:

\begin{Th*}[Théorème 94]
Lorsque le corps cyclique relatif $K$ de degré premier impair $\ell$ par rapport à $k$ a sa différente relative égale à 1, il y a toujours dans $k$ un idéal $I$ qui n'est pas un idéal principal de $k$ mais qui devient un idéal principal dans $K$. La $\ell$-ième puissance de cet idéal $I$ est alors aussi nécessairement un idéal principal dans $k$ et le nombre des classes du corps $k$ est divisible par $\ell$.
\end{Th*}

La démonstration de Hilbert repose sur le fait que, sous les hypothèses énoncées (absence de ramification, groupe de Galois cyclique, degré premier impair) , il existe une unité $\eta$ de $K$, de norme relative $N_{K/k}(\eta) = 1$, qui n'est pas la puissance symbolique $\varepsilon^{\sigma-1}$ d'une autre unité, pour un $\sigma$ non trivial de $G =\Gal (K/k)$ (Théorème 92). Par le Théorème 90, une telle unité s'écrit $\eta = \alpha^{\sigma-1}$ pour un $\alpha$ de $K$, qui engendre ainsi un idéal ambige (i.e. invariant par $G$) de $K$,donc étendu de $k$ puisque $K/k$ est supposée non ramifiée (Théorème 93). Cet idéal ne pouvant être principal dans $k$ (sans quoi on aurait $\alpha=\beta\varepsilon$, pour un $\beta$ de $k$ et une unité $\varepsilon$ de $K$, donc $\eta = \alpha^{\sigma-1}=\varepsilon^{\sigma-1}$ contrairement à l'hypothèse), c'est l'idéal attendu $I$. Naturellement, il vient alors $I^\ell= N_{K/k}(I) = \big(N_{K/k}(\alpha)\big)$, ce qui montre que la classe de $I$ est exactement d'ordre $\ell$.
\smallskip

On voit clairement par là que, dans le résultat de Hilbert, le nœud de la preuve consiste à associer aux idéaux $\a$ de $k$ qui capitulent dans $K$ les applications $f_\a: \sigma\mapsto \alpha^{\sigma-1}$ , où $\alpha$ est un générateur arbitraire de $\a$ dans $K$. En termes modernes, $f_\a$  est ce que nous appelons un 1-cocycle du groupe $G$ à valeurs dans le groupe dans les unités $E_K$. Or, comme l'a remarqué Iwasawa \cite{Iwa}, sous sa forme cohomologique la correspondance obtenue est parfaitement générale.\smallskip

Il vient, en effet:

\begin{Prop} 
Étant donnée une extension galoisienne quelconque $K/$k de corps de nombres, le quotient $P^G_K/P_k$ du groupe des idéaux principaux ambiges de $K$ par le sous-groupe des idéaux principaux de $k$ s'identifie au premier groupe de cohomologie $H^1(G,E_K)$ du groupe de Galois $G = \Gal(K/k)$ à valeurs dans le groupe des unités de $K$.\par

En particulier, lorsque l'extension $K/k$ est non ramifiée, les idéaux ambiges étant étendus, il suit canoniquement:

\centerline{$Cap_{K/k} \simeq H^1(G,E_K)$.}
\end{Prop}

\Preuve Partons de la suite exacte courte qui définit le groupe $P_K$:
$$
1 \longrightarrow E_K \longrightarrow K^\times \longrightarrow P_K \longrightarrow 1.
$$
La suite exacte de cohomologie associée commence ainsi
$$
1 \longrightarrow E_k \longrightarrow k^\times \longrightarrow P_K^G\longrightarrow H^1(G,E_K) \longrightarrow H^1(G,K^\times)
$$
Et le théorème 90 de Hilbert généralisé nous dit que terme de droite $H^1(G,K^\times)$ est nul. II vient donc, comme attendu:
$$
P^G_K/ P_k \simeq  H^1(G,E_K).
$$ 
Cela étant, si l'extension $K/k$ est non ramifiée, un calcul élémentaire (cf. eg. \cite{Che}, p. 403) montre que les idéaux ambiges de $K$ sont exactement les étendus à  $K$ des idéaux de $k$. Le groupe  $P_K^G$ est donc constitué des idéaux de $k$ qui capitulent dans $K$, ce qui donne le résultat annoncé\footnote{Lorsque l'extension $K/k$ se ramifie, l'identité $P_K^G = P_K \cap I_k$ est en défaut, et la capitulation $Cap_{K/k}$ apparaît alors comme noyau de l'application naturelle $H^1(G,E_K) \to I_K^G//I_k$ (cf. \cite{J18}, Ch. III, §1).}.\smallskip

L'isomorphisme obtenu permet de ramener dans tous les cas l'étude de la capitulation à celle de la cohomologie des unités. Le malheur est que celle-ci est extrêmement mal connue, sauf, peut être lorsque le groupe G est cyclique, auquel cas la colomologie des unités est elle même cyclique et l'on a l'identité de Herbrand (cf. \cite{Che}, p. 405):

\begin{Prop}[Quotient de Herbrand des unités]\label{QH}
 Dans une extension cyclique $K/k$ de corps de nombres, le quotient des ordres des groupes de cohomologie galoisienne associés aux unités est donné par la formule :
$$
\frac{|H^0(G , E_K)|}{|H^1(G , E_K)|} =\frac{\prod_{\p_{\si{\infty}}|\si{\infty}}d_{\p_{\si{\infty}}}}{[K:k]}
$$
où, pour chaque place archimédienne $\p_{\si{\infty}}$ du corps $k$, l'entier $d_{\p_{\si{\infty}}}$ désigne le degré de l'extension locale correspondante.
\end{Prop}

Dans la formule obtenue, l'ordre du numérateur $|H^0 ( G , R_K) | = \big(E_k: N_{K/k}(E_K)\big)$ n'est pas connu en général. Cependant,
il vaut au moins 1, ce qui permet d'écrire sous les hypothèses de la proposition 1 :
$$
|Cap_{K/k}|\ge \frac{[K:k]}{\prod_{\p_{\si{\infty}}|\si{\infty}}d_{\p_{\si{\infty}}}}
$$

 Interdisant alors aux places réelles de $k$ de se complexifier dans $K$ (condition qui est automatiquement remplie lorsque le degré $[K:k]$ est impair), nous obtenons le théorème:

\begin{Th}[Théorème 94 généralisé]\label{Th3}
Dans une extension cyclique non ramifiée de corps de nombres, où les places à l'infini sont complètement décomposées, l'ordre de la capitulation est un multiple du degré de l'extension.
\end{Th}

Tout groupe abélien étant produit direct de groupes cycliques, ce résultat suggère la {\em Conjecture principale} suivante:

\begin{Conj}\label{CP}
Dans une extension abéltenne non ramifiée de corps de nombres, où les places à l'infini sont complètement décomposées, l'ordre de la capitulation est un multiple du degré.
\end{Conj}

Nous verrons plus loin, par des arguments tirés du corps de classes, que cette conjecture est vraie lorsque l'extension $K/k$ considérée est maximale, i.e. lorsque $K$ est le corps de classes de Hilbert de $k$. Contentons nous pour l'instant de noter qu'en vertu de l'identité déjà citée
\smallskip

\centerline{$N_{K/k}\circ j_{K/k}=[K/k] $,}\smallskip

\noindent la capitulation dans une extension (abélienne) de degré $n$ n'intéresse que les classes d'ordre divisant $n$; en particulier la restriction de l'homomorphisme d'extension $j_{K/k}$ au sous-groupe de $Cl_k$ formé des classes d'ordre étranger à $n$ est injective. Si donc nous écrivons $G = \prod G_p$ la décomposition du groupe $G=\Gal(K/k)$ comme produit direct de ses sous-groupes de Sylow, et que nous convenons de désigner, pour chaque premier $p$ divisant $n$, par $K^{(p)}$ la p-sous-extension maximale de $K$, la $p$-partie $Cap^{(p)}_{K/k}$ du groupe $Cap_{K/k}$ n'est autre que le sous-groupe $Cap_{K^{(p)}/k}$ qui mesure la capitulation dans la p-extension $K^{(p)}/k$ de groupe $\Gal(K^{(p)}/k)\simeq G_p$.\smallskip

Il suit de là que la conjecture précédente se lit dans les $p$-extensions. Cela étant, nous avons:

\begin{Prop} Soit K/k une extension abélienne non ramifiée et $\infty$-décomposée dont le groupe de Galois $G$ est le produit d'un groupe cyclique par un groupe cyclique élémentaire. Alors l'ordre de la capitulation $Cap_{K/k}$ est un multiple du degré $[K:k]$ de l'extension . 
\end{Prop}

\noindent {\em Preuve} (d'après Adachi \cite{Ada}). D'après ce qui précède, il suffit d'établir la proposition lorsque $K/k$ est une $p$-extension, c'est-à-dire lorsque le groupe de Galois $G = \Gal(K/k)$ est le produit d'un $p$-groupe cyclique par un $p$-groupe cyclique élémentaire $H$. Dans ce cas, la suite exacte longue de HochschiId-Serre, appliquée au groupe des unités de $K$, commence ainsi (avec $L = K^H$):\smallskip

\centerline{$1 \longrightarrow H^1(G/H,E_L) \longrightarrow H^1(G,E_K) \longrightarrow H^1(G,E_K)^H \longrightarrow H^2 (G/H,E_L)$}\smallskip

\noindent II vient donc :

\centerline{$|Cap_{K/k}| = |H^1G,E_K)| \ge |H^1G/H,E_L)| .|H^1(G,E_K)^H|$}\smallskip

\noindent Cela étant, le premier groupe $H^1(G/H,E_L)$ est d'ordre au moins $[L:k]$, d'après le Théorème \ref{Th3}; et le second $H^1(G,E_K)^H$ est au moins d'ordre $p = [K:L]$, puisque c'est le sous-groupe des points fixes d'un $p$-groupe non trivial. Il suit ainsi:\smallskip

\centerline{$|Cap_{K/k}|\ge [K:L][L:k]  = [K:k]$, comme annoncé.}\medskip

En dehors cependant de ce cas particulier, force est de constater que les minorations arithmétiques de la capitulation n'ont guère conduit jusqu'ici à des résultats significatifs, eu égard à l'inégalité espérée en toute généralité. En revanche, comme nous allons le voir, il est possible dans certaines situations spécifiques, d'améliorer la Conjecture \ref{CP}.

\subsection{Le Spiegelungssatz de Leopoldt comme résultat de capitulation}

Revenons sur l'isomorphisme donné par la proposition 1:\smallskip

\centerline{$Cap_{K/k} \simeq H^1(G,E_K)$.}\smallskip

Faute de connaître la cohomologie des unités en toute généralité, il est tentant de regarder ce qui se passe lorsque le groupe $E_K$ se réduit à son sous-groupe de torsion $\mu_K$, ou, moins restrictivement, lorsque $E_K$ contient $\mu_K$ comme facteur direct, en un certain sens de façon canonique. C'est précisément ce qui se produit lorsque le corps $K$ possède une conjugaison complexe naturelle, c'est-à-dire lorsque $K$ est une extension quadratique totalement imaginaire d'un corps de nombres totalement réel $K_+$.\smallskip

Introduisons pour simplifier le complété profini $\hat\ZZ$ de $\ZZ$ pour la topologie définie par ses sous-groupes d'indice impair:\smallskip

\centerline{$\hat\ZZ=\underset{n\;\si{\rm impair}}{\varprojlim}\ZZ/n\ZZ$}\smallskip

(de sorte que $\hat\ZZ$ s'identifie canoniquement au produit $\prod_{p\;\si{\rm impair}}\ZZ_p$
des complétés p-adiques de $\ZZ$ pour tous les nombres premiers impairs).

Si $\Delta=\{1,\tau\}$ désigne le groupe de Galois $\Gal(K/K_+)$, le nombre 2 étant inversible dans $\hat\ZZ$, l'algèbre de Galois $\ZZ[\Delta]$ se décompose comme produit direct de deux exemplaires de $\hat\ZZ$ à l'aide des idempotents orthogonaux
$$
e_+=\frac{1}{2}(1+\tau) \qquad \& \qquad e_-=\frac{1}{2}(1-\tau)
$$
Plus généralement, tout $\hat\ZZ[\Delta]$-module $M$ s'écrit comme somme directe de sa composante {\em réelle} $ M^+ = e_+M$ et de sa composante {\em imaginaire} $M^- = e_-M$. Par exemple, le sous-groupe $\hat{Cl}_K$ de $Cl_K$ formé des classes d'ordre impair s'écrit canoniquement :
$$
\hat{Cl}_K=\hat{Cl}_K^+ \oplus \hat{Cl}_K^-
$$
et le premier facteur $\hat{Cl}_K^+$ n'est autre que le sous-groupe $\hat{Cl}_{K_+}$ des classes d'ordre impair du groupe des classes d'idéaux du sous-corps réel maximal $K_+$. de $K$.\medskip

Supposons maintenant que le groupe $G$ commute à la conjugaison complexe $\tau$, autrement dit que $K$ comme $k$ admettent une conjugaison complexe naturelle, et que l'extension $K/k$ provienne d'une extension galoisienne de même degré $K_+/k_+$ pour leurs sous-corps réels maximaux. Tensorisant par $\hat\ZZ$ les suites exactes utilisées dans la démonstration de la proposition 1, nous obtenons directement l'isomorphisme de $\hat\ZZ[\Delta]$-modules:
$$
\hat C\!ap_{K/k} \simeq H^1(G,\hat E_K)
$$
où $\hat C\!ap_{K/k}$  désigne le sous-groupe de $Cl_K$ formé des classes d'ordre impair qui capitulent dans $K$, et $\hat E_K=\hat\ZZ\otimes_\ZZ E_K \simeq \underset{n\;\si{\rm impair}}{\varprojlim} E_K/E_K^n$ est le complété profini de $E_K$ pour la topologie des sous-groupes d'indice impair.\smallskip

Par action de l'idempotent $e_-$, nous en déduisons immédiatement l'isomorphisme entre composantes imaginaires:
$$
\hat C\!ap{}_{K/k}^- \simeq H^1(G,\hat{E}_K)^- \simeq H^1(G,\hat\mu_K)
$$
puisque la composante imaginaire du groupe $\hat E_K$, qui est $\hat\ZZ$-engendrée par les classes dans $\hat E_K$ des unités de $K$ dont tous les conjugués sont de module 1, se réduit au sous-groupes $\mu_K$ des racines d'ordre impair de l'unité dans $K$.\smallskip

Nous avons même mieux : La démonstration de la proposition 1 utilisait l'hypothèse de non-ramification pour établir l'égalité $Cap_{K/k}\simeq P_K^G/P_k$. Ici, où seules nous intéressent en fin de compte les composantes imaginaires des groupes considérés, il nous suffit que cette égalité ait lieu pour les composantes imaginaires des sous-groupes d'ordre impair, c'est-à-dire en l'occurrence que les idéaux ramifiés dans $K/k$ ne se décomposent pas dans $K/K_+$ . Nous pouvons donc énoncer:

\begin{Th}
Soit $K/k$ une extension galoisienne de corps à conjugaison complexe (autrement dit, soit $K_+/k_+$ une extension galoisienne de corps totalement réels, puis $k$ une extension quadratique totalement imaginaire de $k_+$ et $K_+$ l'extension posée $K_+k$). Si les places ramifiées dans $K/k$ ne sont pas décomposées par la conjugaison complexe (en particulier, si l'extension $K/k$ est non ramifiée), la composante imaginaire du plus grand sous-groupe d'ordre impair de la capitulation $\hat C\!ap_{K/k}$ est donnée par l'isomorphisme
$$
\hat C\!ap_{K/k}^- \simeq H^1(\Gal(K/k),\hat\mu_K),
$$
où $\hat\mu_K$ est le groupe des racines d'ordre impair de l'unité dans $K$.
\end{Th}

\begin{Cor}
Conservons les hypothèses du théorème 6 et sons en outre $\hat\mu_K=\hat\mu_k$  (i.e. que les racines d'ordre impair de l'unité dans $K$ sont déjà dans $k$) . Il vient alors:
$$
\hat C\!ap_{K/k}^- \simeq \Hom (\Gal(K/k),\hat\mu_K),
$$
et la composante imaginaire du plus grand sous-groupe d' ordre impair de la capitulation s' identifie ainsi au radical kummérien de la sous-extension ahélienne maximale d' exposant $n=|\hat\mu_k|$ de $K/k$. 
\end{Cor}

 Dans ce dernier cas, l'ordre de la composante imaginaire de la capitulation coïncide avec le degré de l'extension chaque fois que le groupe de Galois $G=\Gal(K/k)$ est abélien et d'exposant $n$. En particulier l'ordre de la capitulation totale (réelle et imaginaire) est alors strictement supérieur au degré de l'extension. En effet, le théorème 3 appliqué à l'une quelconque des sous-extensions cycliques de $K_+/k_+$ montre que la composante réelle de la capitulation n'est jamais triviale.\smallskip

On voit donc clairement sur cet exemple combien il peut être difficile de fixer en général une borne supérieure à l'ordre de la capitulation dans une extension quelconque (autre que celle résultant de l'inclusion triviale $C\!ap_{K/k} \subset {}_{[K:k]}Cl_k$.\smallskip

Le corollaire 7 ci-dessus peut être regardé comme une généralisation du très classique Spiegelungssatz de Leopoldt (cf. \cite{Leo}). Fixons, en effet, un nombre premier impair $\ell$, et considérons un corps totalement réel $k_\circ$ ; notons $k=k_\circ[\zeta]$ l'extension cyclotomique engendrée sur $k_\circ$ par les racines $\ell$-ièmes de l'unité, puis $k_+=k_\circ[\zeta+\zeta^{\si{-1}}]$ le sous-corps réel maximal de $k$. Écrivons $\ell^d$ (avec $d \ge 1$) l'ordre du $\ell$-groupe des racines de l'unité dans $k$, puis, pour tout $m = 1,\dots,d$, introduisons la $\ell$-extension abélienne non-ramifiée d'exposant $\ell^m$ maximale $K_+^{(m)}$ de $k_+$. Par la théorie du corps de classes, le groupe de Galois $G^{(m)}=\Gal(K_+^{(m)}/k_+)$ s'identifie au quotient d'exposant $\ell^m$ du groupe des classes d'idéaux du corps totalement réel $k_+$.
$$
G^{(m)} \simeq {}^{\ell^{\si{m}}}Cl_{k_+} = Cl_{k_+}/Cl_{k_+}^{\ell^{\si{m}}}
$$
i.e. à la composante réelle du quotient d'exposant $\ell^m$ du groupe des classes du corps à conjugaison complexe $k$:
$$
G^{(m)} \simeq {}^{\ell^{\si{m}}}Cl_k^+.
$$
Maintenant, si $K^{(m)}$ désigne le corps composé $K^{(m)}_+k$, l'extension $K^{(m)}/k$ vérifie les hypothèses du corollaire 7 pour le nombre premier $\ell$ (en ce sens que les $\ell$-sous-groupes de Sylow des groupes$\mu_K$ et $\mu_{K^{\si{(m)}}}$ coïncident), de sorte que nous obtenons:
$$
C\!ap^-{K^{\si{(m)}}/k} \simeq \Hom(G^{(m)},\hat\mu_k) \simeq \Hom(^{\ell^{\si{m}}}Cl_k^+,_{\ell^{\si{m}}}\mu_k).
$$
 Et les classes qui capitulent sont au plus d'ordre $\ell^m$. De l'inclusion évidente $C\!ap^-{K^{\si{(m)}}/k} \subset {}_{\ell^{\si{m}}}\Cl_k^-$, nous déduisons donc finalement l'homomorphisme injectif:
$$
 \Hom(^{\ell^{\si{m}}\!}Cl_k^+,{\;}_{\ell^{\si{m}}}\mu_k) \hookrightarrow {}_{\ell^{\si{m}}}\Cl_k^-
$$
qui nous montre que, pour $m = 1,\dots,d$, le $\ell^m$-rang de la composante réelle du $\ell$-groupe des classes d'idéaux du corps $k$ est inférieur au $\ell^m$-rang de la composante imaginaire; ce qui est une extension classique du résultat de Leopoldt (cf. C34], Ch. I).

\section{L'approche algébrique}

\subsection{Le théorème d'Artin-Furwängler}

La Conjecture \ref{CP} postule que dans une extension abélienne non ramifiée et non-décomposée $K/k$ de corps de nombres, l'ordre de la capitulation $C\!ap_{K/k}$ est un multiple du degré $[K:k]$ de l'extension. Cela doit être le cas en particulier lorsque l'extension $K/k$ est maximale sous les conditions énoncées, c'est-à-dire lorsque $K$ est le corps des classes de Hilbert $H_k$ de $k$. La Théorie du corps de classes nous dit alors que le groupe de Galois $\Gal(K/k)$ s'identifie au groupe des classes d'idéaux $Cl_k$ de $k$. Dans ces conditions, le degré de l'extension est exactement l'ordre de $Cl_k$ et la conjecture 4 affirme simplement que tous les idéaux de $k$ capitulent dans K. C'est le résultat d'Artin-Furtwängler (cf. \cite{Art, Fuw}), qui est connu sous le nom de {\em Théorème de l'idéal principal}:

\begin{Th}[Théorème de l'idéal principal]\label{Th8}
Les idéaux d'un corps de nombres deviennent principaux dans son corps des classes de Hilbert (i.e. dans l'extension abélienne non-ramifiée $\si{\infty}$-décomposée maximale de ce corps).
\end{Th}

La section 2 de la bibliographie est consacrée précisément à ce résultat et à ses extensions les plus immédiates. Les démonstrations "historiques" ont été dans leur presque totalité publiées dans les années trente aux annales de Hambourg. Des exposés plus modernes se trouvent dans les livres de Neukirch \cite{Neu} (qui reprend une argumentation de Witt) et d'Artin-Tate \cite{A-T} (que nous suivrons ici, et qui utilisent un soupçon de cohomologie). Dans tous les cas, la preuve du théorème de l'idéal principal repose sur un argument purement algébrique que nous allons maintenant exposer:\smallskip

Étant donnée une extension $K/k$ de corps de nombres, la Théorie du corps de classes identifie le groupe $Cl_k$ au groupe de Galois $G_k = \Gal(H_k/k)$ de l'extension abélienne non ramifiée $\si{\infty}$-décomposée maximale de $k$, et le groupe $Cl_K$ au groupe de Galois  $G_K = \Gal(H_K/K)$ de l'extension correspondante de $K$, Cela étant, dans le diagramme
\begin{displaymath}
\xymatrix@C=2cm{
\;\Cl_K \ar@{->}[r] ^\sim \ar@<1ex>[d]^{N_{L/K}} & \;\;G_K \ar@<1ex>[d]^{\res_{K/k}} \\
\;\Cl_k \ar@{->}[r] ^\sim \ar@{->}[u]^{j_{L/K}} \ar@{->}[r] ^\sim  & \;\;G_k \ar@{->}[u]^{\Ver_{K/k}}
}
\end{displaymath}
la norme arithmétique $N_{K/k}$ correspond à la restriction $\res_{K/k}$ des automorphismes de Galois, et l'homomorphisme d'extension $J_{K/k}$ au transfert $\Ver_{K/k}$ (Verlagerung).

Si maintenant $K$ est contenu dans $H_k$, nous pouvons résumer la situation par le diagramme:\medskip

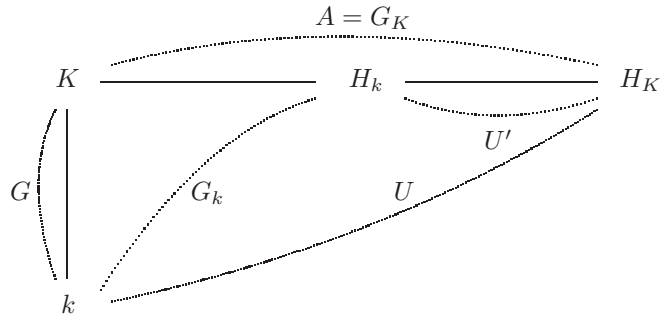
\begin{center}
\unitlength=1.5cm
\begin{picture}(6.6,2.8)

\put(0.65,0){$k$}
\put(0.7,0.3){\line(0,1){1.5}}
\put(0.6,2){$K$}

\bezier{60}(0.6,0.3)(0.3,1.2)(0.6,1.8)
\put(0.2,1.0){$G$}

\put(1.0,2.05){\line(1,0){1.9}}
\put(3.2,2){$H_k$}
\put(3.7,2.05){\line(1,0){1.7}}
\put(5.6,2){$H_K$}

\bezier{120}(1.1,2.2)(3,2.7)(5.4,2.2)
\put(2.9,2.55){$A=G_K$}

\bezier{70}(1.0,0.2)(1.9,1.6)(2.9,1.9)
\put(1.8,1.0){$G_k$}

\bezier{50}(3.7,1.9)(4.5,1.6)(5.4,1.9)
\put(4.4,1.45){$U'$}

\bezier{180}(1.1,0.1)(3.5,0.6)(5.4,1.8)
\put(3.6,1.0){$U$}
\end{picture}
\end{center}
\medskip\medskip

Et, du point de vue purement algébrique, le problème de la capitulation se pose comme suit: étant donnés un groupe $U$ (en l'occurrence $U= \Gal(H_k/k)$) et un sous-groupe abélien $A$ contenant le sous- groupe dérivé $U'$ (en l'occurrence $A = \Gal(H_K/K)$), déterminer le noyau du transfert:\smallskip

\centerline{$\Ver_{U/A}:\; U/U' \to A$.}\medskip

Le Théorème de l'idéal principal résulte de la constatation, purement algébrique, que $\Ver_{U/A}$ est nul lorsque $A$ coïncide avec $U'$.\smallskip

Pour établir ce point, remarquons d'abord que le groupe $U$ se présente comme extension de son sous-groupe normal abélien $A$ par le quotient abélien $G = U/A$:
$$
1 \longrightarrow A \longrightarrow U \longrightarrow G \longrightarrow 1.
$$
À isomorphisme près, la loi de groupe sur $U$ est donc déterminée par la donnée des lois de groupe de $A$ et de $G$, de l'action de $G$ sur $A$ (via les automorphismes intérieurs de $U$) et de la classe dans le groupe de cohomologie $H^2(G,A)$ du système de facteurs 
$$
a_{\sigma,\tau}=u_\sigma u_\tau u_{\sigma\tau}^{-1}
$$
associé au choix d'un système de représentants dans $U$ des classes du quotient $G = U/A$.
Avec ces notations, le groupe $U$ peut être décrit comme l'ensemble des symboles $a.u_\sigma$ (avec $a\in A$ et $\sigma\in G$) équipé de la loi de groupe:
$$
( a . u_\sigma )\,(b.u_\sigma) = ab^\sigma .\, a_{\sigma,\tau}  u_{\sigma\tau} .
$$
Cela étant, le transfert $\Ver_{U/A}$ de $U/U'$ vers $A$ est défini par la formule (cf. \cite{A-T}, Ch. XIII, §2):
$$
\Ver_{U/A} ( a . u_\sigma U') = \prod_{\sigma\in G} \big(u_\sigma(a.u_\tau)u_{\sigma\tau}^{-1}\big) = \prod_{\sigma\in G} a^\sigma a_{\sigma,\tau}
$$
ce que nous pouvons encore écrire:

$$
\Ver_{U/A} ( a . u_\sigma U') = N(a) \prod_{\sigma\in G} a_{\sigma,\tau}
$$
en notant $N(a) = \prod_{\sigma\in G} a^\sigma$ le produit des conjugués de $a$.\smallskip

 Pour réinterpréter le transfert en termes d'algèbre linéaire, introduisons le module résolvant
$$
B= A \oplus I_G = A \oplus \big(\underset{\tau\ne1}{\oplus}\ZZ\,(\tau-1)\big)
$$
somme directe du groupe abélien $A$, regardé comme $\ZZ$-module et noté additivement, et de l'idéal d'augmentation $I_G$ de l'algèbre $\ZZ[G]$; et munissons le de l'action de $G$ définie par :
$$
\sigma * a=a^\sigma \qquad \& \qquad \sigma *(\tau-1)= a_{\sigma,\tau} + \sigma(\tau-1).
$$

Un calcul élémentaire montre alors que l'application $\log : a . u_\tau U' \mapsto a+(\tau-1) + I_G*B$ est un isomorphisme du groupe $U$ abélianisé sur le quotient $B/I_G*B$:
$$
U / U' \simeq B/IG*B .
$$
Cela étant, nous avons:

\begin{LemA}
Dans l'isomorphisme précédent, le transfert $Ver_{U/A}$ correspond à l'application trace:\smallskip

\centerline{ $\Tr_{B/A}=\sum_{\sigma\in G}\sigma$.}
\end{LemA}

\Preuve Il s'agit de vérifier la commutativité du diagramme (où les notations sont multiplicatives à gauche et additives à droite):

\begin{displaymath}
\xymatrix@C=1.5cm{
U/U' \ar@{->}[r] ^\log \ar@{->}[d]_{\Ver_{U/A}} & B/I_G *B \ar@{->}[d]^{\Tr_{B/A}}\\
A \ar@{=}[r] & A
}
\end{displaymath}

Or, nous avons directement, en notations multiplicatives dans $U$:
$$
\Ver_{U/A} ( a . u_\sigma U') = \prod_{\sigma\in G} a^\sigma  \prod_{\sigma\in G} a_{\sigma,\tau}
$$
c'est-à-dire, en notations additives dans $B$:
$$
\Ver_{U/A} ( a . u_\sigma U') =  \sum_{\sigma\in G} a^\sigma \; +\; \sum_{\sigma\in G} a_{\sigma,\tau} 
= \sum_{\sigma\in G} \sigma*a +\sum_{\sigma\in G}\sigma*(\tau-1)\;-\;\big(\sum_{\sigma\in G}\sigma\big)*(\tau-1)
$$
et finalement:
$$
\Ver_{U/A} ( a . u_\sigma U') = \big(\sum_{\sigma\in G}\sigma\big)*\big(a+(\tau-1)\big) = \Tr_{B/A}\big(a+(\tau-1)\big),
$$
comme attendu.\medskip

Et le théorème d'Artin-Furtwängler résulte alors de l'identité $(A:U') \Tr_{B/A}(B) = 0$, qui montre que la trace est nulle sur $B$ pour $A=U'$ et que nous allons maintenant établir:

\begin{LemB}Avec les notations précédentes, on a: $(A:U') \Tr_{B/A}(B) = 0$.
\end{LemB}

\Preuve Partons de la définition du module résolvant:
$$
B=A\oplus I_G
$$
L'isomorphisme $B/I_G*B \simeq U/U'$ montre que le sous-module de torsion $A \cap I_G*B$ de $I_G*B$ est égal à $U'$, ce qui permet d'écrire (pour la structure de Z-module):
$$
I_G*B = U' \oplus I_G^2.
$$
En particulier, il vient :
$$
I_G/I_G^2 \simeq B/(A+I_G*B) \simeq U/A = G .
$$

Décomposons donc le groupe $G$ comme produit direct de sous-groupes cycliques $G_i$ (avec, disons, $i=1,\dots,s$) d'ordres respectifs $e_i$; et faisons choix pour chaque indice $i$ d'un générateur $\tau_i$ dans $G_i$. L'élément $\tau_i$ est ainsi d'ordre $e_i$ dans $G$, et il en est donc de même de l'élément $b_i = \tau_i-1$ dans $B/I_G*B$. Complétons maintenant les $b_i$ déjà obtenus en choisissant des générateurs $b_{s+1},\dots,b_m$ du groupe $A$ d'ordres respectifs $e_i$ modulo $U'$ de telle sorte que nous ayons:\smallskip

\centerline{$\prod_{i=s+1}^m e_i = (A:U')$ et, bien entendu: $\prod_{i=1}^s e_i = (B:A+I_G*B)=(U:A)$.}\smallskip

Par construction, nous avons donc $e_i b_i \in I_G*B, \quad \forall i=1,\dots,m$; ce que nous pouvons écrire:
$$
e_i b_i = \sum_{j=1}^m \theta_{ij}*b_j
$$
pour des $\theta_{ij}$ convenables dans $I_G$, puisque les $b_i$ dans leur ensemble forment un système générateur du $\ZZ[G]$-module $B$. Introduisons alors la matrice $M$ de terme générique:\smallskip

\centerline{$m_{i j} =e_i\delta_{ij}$ (où $\delta_{ij}$ est le symbole de Kronecker).}\smallskip

Nous pouvons réécrire l'identité précédente sous la forme: $\sum_{j=1}^m m_{ij}*b_j=0$; ce qui, par multiplication à gauche par la transcomatrice $\wi M$ de $M$, nous donne:\smallskip

\centerline{$\det M *b_j=0$ pour $j=1,\dots,m$, c'est-à-dire, finalement:}\smallskip

\centerline{$\det M * B =0$.}\smallskip

Reste à évaluer l'élément $d = \det M$ de l'algèbre $\ZZ[G]$. L'isomorphisme $B/A \simeq I_G$ montre que les seuls éléments $d$ de $\ZZ[G]$ qui vérifient $d*B\subset A$ sont les multiples de la trace $\Tr_{B/A}=\sum_{\sigma\in G}\sigma$. Tout le problème est donc d'estimer le degré de $d$. Mais celui-ci est, tout simplement, le déterminant de la matrice en les degrés des $m_{ij}$. Et l'identité $\deg m_{ij}= \delta_{ij} e_i$ donne donc directement:
$$
\deg d = \prod_{i=1}^m e_i = \prod_{i=1}^s e_i \prod_{i=s+1}^m e_i = (U:A)(A:U') = (A:U') \deg\Tr_{B/A}.
$$
Il vient donc: $d=(A:U')\Tr_{B/A}$; ce qui achève la démonstration du Théorème \ref{Th8}.

\subsection{Le théorème de Tannaka-Terada}

Le théorème de Tannaka-Terada est une généralisation de celui d'Artin-Furtwängler qui fait intervenir la théorie des genres, mais repose comme celui-ci sur un résultat purement algébrique de théorie des groupes. Il a été pressenti par Tannaka qui l'a établi d'abord dans un certain nombre de cas particuliers avant que Terada n'en donne finalement une preuve générale quoique extraordinairement technique (cf. \cite{T-T,Te1,Ta2,Ta3,Ta4}). Nous nous inspirons ici de la démonstration extrêmement astucieuse mise au point plus tard en collaboration avec  Adachi (cf. \cite{Ada}) pour établir directement une version sensiblement plus forte du résultat algébrique qui sous-tend le théorème arithmétique, et que Miyake a récemment déduite du résultat classique de Terada (cf. \cite{Mi5}).\smallskip

Commençons par énoncer le résultat historique de Tannaka-Terada :

\begin{Th}
Pour toute sous-extension cyclique $K$ du corps des classes de Hilbert $H_k$ d'un corps de nombres $k$, les classes d'idéaux du corps $K$ qui sont ambiges pour l'action de $\Gal(K/k)$ capitulent dans $H_k$.
\end{Th}

Bien entendu, prenant $K = k$, on retrouve le théorème d'Artin-Furtwängler.\smallskip

{\em Interprétation algébrique.} Notons $L = H_k$ le corps de classes de Hilbert de $k$, puis $H_K$ celui de $K$, et $H_L$ celui de $L$. Écrivons $\Gamma$ le groupe de Galois $\Gal(K/k)$, et faisons choix d'un générateur $\gamma$ de $\Gamma$, ou plutôt d'un relèvement de ce générateur dans $\Gal(H_L/k)$. La théorie des genres (cf. \cite{J18}, Ch. III, §2) nous enseigne que, puisque l'extension $K/k$ est cyclique, la sous-extension maximale $L = H_k$ de $H_K$ qui est abélienne sur $k$ (autrement dit le corps des genres de $H_K$ relativement à $K/k$) est caractérisée comme extension de $K$ par l'identité:
$$
\Gal(L/K) = G_K/G_K^{\gamma-1}
$$
où $G_K/G_K^{\gamma-1}$ est le plus grand quotient de $G_K = \Gal(H_K/K)$ sur lequel $\Gamma$ opère trivialement. La situation peut donc se résumer par le diagramme, avec $A/U' = G_K^{\gamma-1}= (U/U')^{\gamma-1}$:\medskip

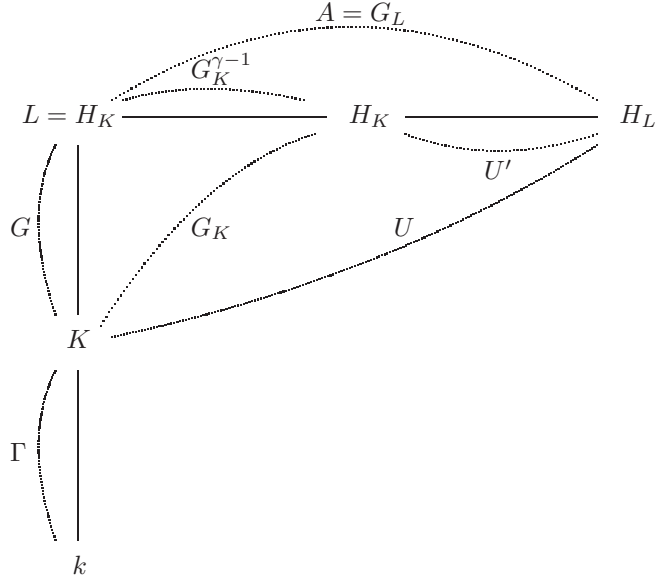
\begin{center}
\unitlength=1.5cm
\begin{picture}(6.6,5.1)

\put(0.75,0){$k$}
\put(0.8,0.3){\line(0,1){1.5}}

\put(0.7,2){$K$}
\put(0.8,2.3){\line(0,1){1.5}}
\put(0.3,4){$L=H_K$}

\bezier{60}(0.6,2.3)(0.3,3.2)(0.6,3.8)
\put(0.2,3.0){$G$}

\bezier{60}(0.6,0.3)(0.3,1.2)(0.6,1.8)
\put(0.2,1.0){$\Gamma$}

\put(1.2,4.05){\line(1,0){1.8}}
\put(3.2,4){$H_K$}
\put(3.7,4.05){\line(1,0){1.7}}
\put(5.6,4){$H_L$}

\bezier{120}(1.1,4.2)(3.2,5.5)(5.4,4.2)
\put(2.9,4.9){$A=G_L$}

\bezier{70}(1.0,2.2)(1.9,3.6)(2.9,3.9)
\put(1.8,3.0){$G_K$}

\bezier{50}(3.7,3.9)(4.5,3.6)(5.4,3.9)
\put(4.4,3.5){$U'$}

\bezier{50}(1.2,4.2)(1.9,4.4)(2.8,4.2)
\put(1.8,4.4){$G_K^{\gamma-1}$}

\bezier{180}(1.1,2.1)(3.5,2.6)(5.4,3.8)
\put(3.6,3.0){$U$}
\end{picture}
\end{center}
\medskip\medskip

Et, du point de vue algébrique, le problème de la capitulation se pose comme suit: étant donnés un groupe $U$ (en l'occurrence $U=\Gal(H_K/K)$), un sous-groupe normal abélien $A$ contenant le sous-groupe dérivé $U'$ (en l'occurrence $A = \Gal(H_L/L)$), et un automorphisme $\gamma$ de $U$ (en l'occurrence, l'automorphisme intérieur associé à un relèvement dans $U$ d'un générateur du groupe cyclique $\Gamma = \Gal(K/k)$), il s'agit de montrer que le noyau du transfert contient le sous- groupe $(U/U')^\Gamma$ des points fixes de $\Gamma$, sous la condition\smallskip

\centerline{$A/U' = (U/U')^{\gamma-1}$, i.e. $A=U^{\gamma-1}U'$.}\smallskip

\noindent Nous allons voir que ce résultat vaut lors même que $\gamma$ n' est pas un automorphisme du groupe $U$:

\begin{Th}\label{Th10}
Soit $A$ un sous-groupe normal abélien d'un groupe $U$, contenant le sous-groupe dérivé $U'$. S' il existe un endomorphisme $\gamma \in   \End U$ tel qu'on ait:\smallskip

\centerline{$A=U^{\gamma-1}U'$,}\smallskip

\noindent alors le transfert $\Ver_{U/A}: U/U'\to A$ est nul sur le sous-groupe $(U/U')^{\langle\gamma\rangle}$ des points fixes de $\gamma$ dans $U/U'$.
\end{Th}

\noindent{\em Reformulation en termes d'algèbre linéaire.} Reprenons la description du groupe $U$ comme extension du sous-groupe abélien $A$ par le quotient abélien $G = U/A$. Le sous-groupe dérivé $U'$ étant évidemment  stable pour l'action de $\gamma$, l'hypothèse$A=U^{\gamma-1}U'$ montre que l'endomorphisme $\gamma$ agit trivialement sur le quotient $U/A = G$, de sorte que nous pouvons décrire l'opération de $\gamma$ sur $U$ en écrivant:\smallskip

\centerline{$u_\tau^\gamma=a_\tau.u_\tau$ pour tout $\tau$ de $G$,}\smallskip

\noindent avec $a_\tau$ convenable dans $A$; ce qui nous permet de transporter cette action au module résolvant $B=A \oplus I_G$ introduit plus haut, en posant\smallskip

\centerline{$\gamma * (\tau-1) = a_\tau + (\tau-1)$.}\smallskip

Notons, pour simplifier $\Gamma = \gamma^\ZZ$ le groupe cyclique engendré par $\gamma$, et $\delta =\gamma-1$. Cela étant, d'après l'isomorphisme $U/U' \simeq B/I_G*B$, l'identité $(U/U')^\delta = A/U'$ s'écrit dans $B$:\smallskip

\centerline{$A + I_G*B = \delta *B + I_G^B$.}\smallskip

\noindent Et tout le problème est de vérifier que l'application trace $\Tr_{B/A} = \sum_{\sigma\in G}\sigma$ de $B/I_G*B$ dans $A$ s'annule sur le sous-module $( B / I_G*B )^\Gamma$ des points fixes par $\Gamma$, autrement dit que l'on a:
$$
\delta b \in I_G*B \Rightarrow \Tr_{B/A}(b) = 0,
$$
conformément au diagramme commutatif (où les notations sont multiplicatives à gauche et additives à droite):

\begin{displaymath}
\xymatrix@C=1.5cm{
(U/U')^\Gamma \ar@{->}[r] ^\sim \ar@{->}[d]_{\Ver_{U/A}} & (B/I_G *B)^\Gamma \ar@{->}[d]^{\Tr_{B/A}}\\
A \ar@{=}[r] & A
}
\end{displaymath}
\smallskip

\noindent{\em Réduction à une identité algébrique.} Les hypothèses faites impliquent que le groupe $U$ est nilpotent puisque son sous-groupe dérivé est abélien. En particulier il est produit direct de ses sous- groupes de Sylow et ce n'est donc pas restreindre la généralité que de supposer que tous les groupes considérés sont des $p$-groupes. Cela étant, nous avons construit dans la section précédente un système de générateurs $b_1,\dots,b_m$ du $\ZZ[G]$-module $B$, tels que les $s$ premiers $b_1,\dots,b_s$ forment une pseudo-base du quotient $B/(A+I_G*B)\simeq I_G/I_G^2$, en en ce sens que leurs ordres respectifs  $e_1,\dots,e_s$ modulo $A+I_G*B$ vérifient la formule:\smallskip

\centerline{$\prod_{i=1}^se_i= \big(B:(A+I_G*B)\big)= (U:A:)$}\smallskip

Si $\Gamma$ est lui-même un $p$-groupe, l'hypothèse $A + I_G*B = \delta *B + I_G^B$ nous assure ici que $b_1,\dots,b_s$ engendrent $B$ comme $\ZZ[G\times\Gamma]$-module. Sinon, quitte à compléter les $b_i$ déjà obtenus par des éléments de $A$, nous pouvons toujours nous ramener à ce cas sans modifier le produit des $e_i$. En particulier l'identité $e_ib_i\in A+I_G*B$ peut donc s'écrire pour chaque $i=1,\dots,s$:\smallskip

\centerline{$e_ib_i=\sum_{j=1}^s \theta_{ij}*b_j + \delta *\sum_{j=1}^s \lambda_{ij}*b_j$ avec $\theta_{ij}\in I_G$ et $\lambda_{ij}\in \ZZ[G\times\Gamma]$.}\smallskip

\noindent Posons, comme précédemment,\smallskip

\centerline{$m_{i j} =e_i\delta_{ij}-\theta_{ij}$ (où $\delta_{ij}$ est le symbole de Kronecker).}\smallskip

\noindent Notons enfin $M$ la matrice $[m_{i j}]$ et $\Lambda$ la matrice $[\lambda_{i j}]$. Nous obtenons:

\begin{LemC}
L'opérateur $\Delta = \det (M-\delta\Lambda) = \Tr_{B/A}-\delta D$ annule le module $B$. En particulier, comme opérateur sur $B$, l'application trace se factorise par $\delta$; ce que nous écrivons:\smallskip

\centerline{$\Tr_{B/A}=\delta D$.}
\end{LemC}

\Preuve Nous avons par construction:\smallskip

\centerline{$\sum_{j=1}^s(m_{ij}-\delta\lambda_{ij})b_j=0$ pour $i=1,\dots,s$; autrement dit: $(M-\delta\Lambda)*{\mathbf b}$,}\smallskip

\noindent où $\mathbf b$ désigne le vecteur colonne formé des $b_i$ pour $i=1,\dots,s$. Multipliant à gauche cette identité par la transcomatrice de $M-\delta\Lambda$, nous obtenons, comme annoncé:\smallskip

\centerline{$\Delta *b_i=0$ pour $ i=1,\dots,s$, c'est-à-dire, en fin de compte: $\Delta * B=0$,}\smallskip

\noindent puisque les $b_i$ engendrent $B$. Maintenant, le développement de $\Delta$ suivant les puissances de $\delta$ est de la forme:

\centerline{$\delta = \det M - \delta D$}\smallskip

\noindent avec $D = \sum_{j=1}^s \det (M_1, \dots,M_{j-1},\Lambda_j, \dots,M_{j+1},\dots,M_s-\delta\Lambda_s)$, si $M_1,\dots,M_s$ d'une part et $\Lambda_1,\dots,\Lambda_s$ d'autre part désignent les vecteurs colonnes respectifs de $M$ et de $A$.

Et tout revient à vérifier que $\det M$ est l'opérateur trace $\Tr_{B/A} = \sum_{\sigma\in G}\sigma$. Or, cela est clair puisque l'identité$(M-\delta\Lambda)*\mathbf b$ lue modulo $A$ s'écrit dans $I_G\simeq B/A$:\smallskip

\centerline{$M *\mathbf  b = 0$ (puisque les $\delta * b_i$  sont dans $A$),}\smallskip

\noindent ce qui donne (après multiplication par la transcomatrice de $M$):\smallskip

\centerline{$\det M * \mathbf  b =0$, i.e. $\det M *B\subset A$}\smallskip

\noindent Par un argument déjà utilisé, $\det M$ qui est dans $\ZZ[G]$ et annule $B/A = I_G$ est donc un multiple de la trace ; et il est égal à la trace puisque nous avons trivialement:\smallskip

\centerline{$\deg (\det M) = \det [\deg m_{ij}] = \prod_{i=1}^s e_i= \big(B: (A+I_G*B) \big) = (U:A)= \deg\Tr_{B/A}$.}\smallskip

En conclusion, pour établir que la trace est nulle sur le sous-groupe $(B/I_G*B)^\Gamma$ de$B/I_G*B$ ou, si l'on préfère, sur le sous-module ${}^{\delta^{\si{-1}}\!}(I_G*B)$ de $B$, il convient de montrer que l'opérateur $D$ est nul sur le sous-module $\delta\big({}^{\delta^{\si{-1}}\!}(I_G*B)\big)= (I_G*B)\cap A = U'$. Mais comme le sous-groupe dérivé $U'$ est engendré par les commutateurs $[a,u_\tau]=a^{\tau-1}$, qui sont dans $I_G*A$, et ceux construits sur les représentants $[u_\sigma,u_\tau] = u_\sigma u_\tau u_\sigma^{-1}u_\tau^{-1}=a_{\sigma,\tau}-a_{\tau,\sigma}$, que l'on peut engendrer sur $\ZZ[G]$ à partir des éléments $[b_h,b_k]=b_h*b_k-b_k*b_h$ , tout le problème se ramène à vérifier l'identité:\smallskip

\centerline{$Db_h*b_k=Db_k*b_h$}\smallskip

\noindent Lorsque celle-ci a lieu, l'hypothèse de départ $A=\delta*B+U'$ nous donne, en effet, immédiatement:\smallskip

\centerline{$D*U'=D I_G*A=D I_G(\delta *B+U')=\Tr_{B/A}(I_G*B)+I_G D*U'=I_G*(D*U')$;}\smallskip

\noindent donc $D*U'=0$, de sorte que le Théorème \ref{Th10} résulte finalement du lemme:

\begin{LemD}
Avec les notations précédentes il vient: $Db_k* b_h = Db_h*b_k$, pour $h \ne k$.
\end{LemD}

\Preuve Partons de l'identité $\sum_{j=1}^s(m_{ij}-\delta\lambda_{ij})b_j=0$ pour $i=1,\dots,s$ dans $B$; et regardons la modulo $A$. Nous obtenons l'identité:\smallskip

\centerline{$\sum_{j=1}^s m_{ij}.b_j=0$ pour $i=1,\dots,s$ dans $I_G=B/A$.}\smallskip

\noindent Cela étant, par un calcul purement formel dans l'algèbre des polynômes $I_G[\delta]$ en l'indéterminée $\delta$, nous avons, en introduisant les vecteurs colonnes:\smallskip

\centerline{$\delta Db_k=\Tr_{B/A}(b_k)-\Delta b_k=-\det (M-\delta\Lambda)b_k=-\det (M_1-\delta\Lambda_1,\dots,(M_k-\delta\Lambda_kb_k,\dots,M_s-\delta\Lambda_s)$,}\smallskip

\noindent c'est-à-dire:\smallskip

\centerline{$\delta Db_k=-\det (M_1-\delta\Lambda_1,\dots,\sum_{j=1}^s M_jb_j-\delta,\sum_{j=1}^s \Lambda_jb_j,\dots,M_s-\delta\Lambda_s)$,}\smallskip

\noindent avec $\sum_{j=1}^s M_jb_j=0$, comme indiqué ci-dessus, d'où:\smallskip

\centerline{$\delta Db_k=\det (M_1-\delta\Lambda_1,\dots,\delta,\sum_{j=1}^s \Lambda_jb_j,\dots,M_s-\delta\Lambda_s)$,}

\noindent et finalement:

\centerline{$Db_k=D_k$ avec $D_k=\det (M_1-\delta\Lambda_1,\dots,\delta,\sum_{j=1}^s \Lambda_jb_j,\dots,M_s-\delta\Lambda_s)$.}\smallskip

\noindent Considérons donc la matrice $N_k= [M_1-\delta\Lambda_1,\dots,\delta,\sum_{j=1}^s \Lambda_jb_j,\dots,M_s-\delta\Lambda_s]$ de déterminant $D_k$ et notons $\wi N_k$ sa transcomatrice.

L'identité de départ, $\sum_{j=1}^s(m_{ij}-\delta\lambda_{ij})b_j=0$ pour $i=1,\dots,s$, s'écrit encore:
$$
N_k *\left[\begin{array}{c}b_1\\ \vdots  \\0 \\ \vdots \\ b_s \end{array}\right] = -(M_k-\delta\Lambda_k)b_k \;;
$$
ce qui, après multiplication à gauche par la transcomatrice $\wi N_k$ de $N_k$, donne:
$$
\left[\begin{array}{c}D_k*b_1\\ \vdots  \\0 \\ \vdots \\ D_k*b_s \end{array}\right]  = 
\wi N_k N_k *\left[\begin{array}{c}b_1\\ \vdots  \\0 \\ \vdots \\ b_s \end{array}\right]  =
-\wi N _k(M_k-\delta\Lambda_k)*b_k .
$$
Ainsi, pour $h \ne k$, le terme  $D_k*b_h$ est obtenu en multipliant la $h$-ième ligne de l'opposé de la matrice $\wi N_k$ par la matrice $M_k-\delta\Lambda_k$, et en faisant agir tout sur $b_k$. Mais la $h$-ième ligne de $\wi N_k$ se calcule en formant les mineurs de la matrice (où l'on a isolé les colonnes d'indices $k$ et $h$):\smallskip

\centerline{$ [M_1-\delta\Lambda_1,\dots,\sum_{j=1}^s \Lambda_jb_j,\dots,\,.\,,\dots,M_s-\delta\Lambda_s]$.}\smallskip

\noindent Elle correspond donc à la $k$-ième ligne de $\wi N_h$, qui se calcule, elle, en formant les mineurs de la matrice:

\centerline{$ [M_1-\delta\Lambda_1,\dots,\,.\,,\dots,\sum_{j=1}^s \Lambda_jb_j,,\dots,M_s-\delta\Lambda_s]$.}\medskip

La multiplication par la matrice colonne $M_k-\delta\Lambda_k$ redonne donc exactement le déterminant $D_k$, d'où comme annoncé:
$$
D_k*b_h=D_h*b_k
$$
ce qui achève la démonstration.

\begin{Cor}\label{C11}
Sous les hypothèses du Théorème, l'ordre du noyau du transfert $\Ver_{U/A}$ est un multiple de l' indice $(U:A)$ du sous-groupe $A$.
\end{Cor}

\Preuve Le Théorème \ref{Th10} affirme en effet que le transfert $\Ver_{U/A}$ est nul sur le sous-groupe $(U/U')^\Gamma$ des points fixes de  dans U/U' Or l'ordre de ce groupe est égal à l'indice de A dans U en vertu de l'exactitude de la suite longue de groupes finis:
$$
\xymatrix{1 \ar@{->}[r] & (U/U')^\Gamma \ar@{->}[r] & U/U' \ar@{->}[r]^{\gamma-1\;} & U/U' \ar@{->}[r] & U/A \ar@{->}[r] & 1}
$$

\section{Interprétation arithmétique des résultats algébriques}

\subsection{Le théorème de l'idéal principal}

Revenons d'abord sur le théorème d'Artin-Furtwängler. Comme l'a observé Herbrand \cite{Her}, le caractère purement formel du Théorème de l'idéal principal permet de le transposer facilement {\em mutatis mutandis} dans le cadre plus vaste des groupes de classes de rayons. On peut même, à l'instar de Miyake \cite{Mi4} l'énoncer de façon très générale en termes idèliques. Introduisons  pour cela le formalisme profini de la Théorie du corps de classes développé dans \cite{J18} sous sa forme $\ell$-adique : Un corps de nombres $K$ étant donné, pour chaque place $\p$ de $K$ notons $\pi_\p$ une uniformisante locale puis $\R_{K_\p}=\varprojlim K_\p^\times/K_\p^{\times\,n}\simeq \U_{K_\p}\,\pi^{\hat \ZZ}$ le compactifié profini du groupe multiplicatif du complété $\p$-adique de $K$, et appelons groupe des idèles de $K$ le produit restreint
$$
\J_K=\prod_\p^{\si{\rm res}}\R_{K_\p} =\bigcup_{S {\;\si{\rm fini}}}\J_K^S\quad  {\rm avec}\quad \J_K^S=\prod_{\p\in S}\R_{K_\p}\prod_{\p\notin S}\U_{K_\p}
$$
constitué des familles d'éléments des $\R_{K_\p}$ dont presque tous sont dans le sous-groupe unité $\U_{K_\p}$, équipé de sa topologie naturelle de limite inductive des sous-groupes compacts $\J_K^S$. Notons enfin
$$
\R_K=\hat\ZZ\otimes_\ZZ K^\times
$$
le tensorisé multiplicatif du groupe $K^\times$ par le groupe profini $\hat\ZZ=\varprojlim \ZZ/n\ZZ$. Dans l'isomorphisme du corps de classes, le quotient
$$
\C_K=\J_K/\R_K
$$
s'identifie (algébriquement et topologiquement) au groupe de Galois $\G_K=\Gal(K^{\si{\rm ab}}/K)$ de l'extension abélienne maximale de $K$.

Cela posé, l'homomorphisme d'extension $j_{K/k}$ pour les groupes d'idèles satisfait la propriété:

\begin{Th}\label{Th12}
 Soit $K/$k une extension galoisienne de corps de nombres de groupe de Galois $G$. Alors, si $\V_K$ est un sous-groupe de $\J_K$ qui contient $\R_K$ et est stable pour l'action de $G$, on a:
$$
\V_K{}^{\si{-1}\!}N_{K/k}(\R_K)=\J_K \quad\Rightarrow\quad j_{K/k}(\J_K)\subset\V_K
$$
\end{Th}

\noindent {\em Preuve} (d'après Miyake \cite{Mi1}). Supposons $\V_K{}^{\si{-1}\!}N_{K/k}(\R_K)=\J_K$ et introduisons l'extension abélienne $L$ de $K$ associée à $\V_K$ par la théorie du corps de classes. Le sous-groupe de $\J_k$ qui correspond à la sous-extension maximale $F$ de $L$ qui est abélienne sur $k$ étant l'image de la norme $N_{K/k}(\V_K)$ dans $\C_k=\J_K/\R_k$, l'identité dans $\J_k$:
$$
\R_k \,N_{K/k}(V_K) = \R_k \,N_{L/k}(\J_K) 
$$
donnée par l'hypothèse nous assure que $F$ est encore la plus grande sous-extension de $K$ abélienne sur $k$. Maintenant, le groupe $\V_K$ étant supposé stable pour l'action de $G = \Gal(K/k)$, le corps $L$ est une extension galoisienne de $k$, et nous pouvons appliquer le théorème de Furtwängler au groupe $U = \Gal(L/k)$ de groupe dérivé $U' = \Gal(L/F)$.

Or, par le corps de classes, nous avons $U/U'\simeq \J_k/\R_k N_{K/k}(\J_K)$ et $U/U''\simeq \J_F/\R_f N_{K/F}(\V_K)$. Et le transfert $\Ver{U/U'}$ correspond à l'extension des idèles de $k$ à $F$. Il vient donc, comme annoncé:
$$
j_{K/k}(\J_k) =j_{K/F}\big( j_{F/k}(\J_k)\big)=j_{K/F}\big(\R_F N_{K/F}((\V_k)\big).
$$

\begin{Sco}
Conservons les notations du Théorème et écrivons $\big(\J_K:{}^{\si{-1}\!}N_{K/k}(\R_k)\V_K\big)=d$. Cela posé, nous obtenons simplement $j_{K/k}(\J_K)^d\subset\V_K$ en dehors de toute hypothèse sur $d$.
\end{Sco}

\Preuve Soit $L$ l'extension abélienne de $K$ correspondant au sous-groupe ouvert ${}^{\si{-1}}N_{K/k}(\R_k)\V_K$. Par construction, nous avons donc $d = [L:K]$ et $\R_KN_{L/K}(\V_K)={}^{\si{-1}}N_{K/k}(\R_k)\V_K$. Nous pouvons dès lors appliquer le Théorème 11 au sous-groupe ouvert ${}^{\si{-1}}N_{K/k}(\V_k)\subset \J_K$ dans l'extension galoisienne $L/k$.
Nous obtenons ainsi les inclusions:
$$
j_{L/k}(\J_k) \subset {}^{\si{-1}\!}N_{K/k}(\V_k),\qquad {\rm puis}\qquad j_{L/k}(\J_k^d)=N_{L/K}\circ j_{L/k}(\J_k) \subset \V_K,
$$
comme attendu.\medskip

\noindent Appliqué maintenant aux groupes de classes de rayons, le Théorème \ref{Th12} peut s'énoncer comme suit:

\begin{Cor}\label{C14}
Soient $K$ un corps de nombres algébriques et $\f_K$ un diviseur entier de $K$. Notons $L = H_K^{\f_{\si{K}}}$ le corps des classes de rayons modulo $\f_K$ sur $K$ (de sorte que le groupe de Galois $\Gal(L/K)$ s'identifie au quotient $Cl_K^{\f_{\si{K}}} = D_K^{\f_{\si{K}}}/R_K^{\f_{\si{K}}}$ du groupe des diviseurs de $K$ étrangers à $\f_K$ par le sous-groupe des diviseurs principaux engendrés par les éléments congrus à 1 modulo $\f_K$) et définissons un diviseur $\f_L$ de $L$ en posant pour chaque place $\p_L$ de $L$ au-dessus d'une place $\p_K$ de $K$ représentée dans $\f_K$:

\centerline{$v_{\p_{\si{L}}}(\f_L)=\Psi_{\p_{\si{L}}}(v_{\p_{\si{K}}}(\f_K)) $,}\medskip

\noindent où $\Psi_{\p_{\si{L}}}$ désigne la fonction de Herbrand associée à la localisée en $\p_L$ de l' extension $L/K$.\par

Cela étant, l'application naturelle $Cl_K^{\f_{\si{K}}} \to Cl_L^{\f_{\si{L}}}$ est l'application nulle. Autrement dit, les classes de rayons de  $Cl_K^{\f_{\si{K}}}$ capitulent dans le groupe $Cl_L^{\f_{\si{L}}}$ des classes de rayons de $L$ modulo $\f_L$.
\end{Cor}

Bien entendu, prenant $\f_K=1$, on obtient $\f_L=1$, ce qui redonne le Théorème \ref{Th8}.\smallskip

\Preuve Rappelons que le corps des classes de rayons de conducteur $\f_K$ est la plus grande extension abélienne de $K$ de conducteur $\f_K$, et qu'elle est associée au groupe d'idèles:\smallskip

\centerline{$\R_KN_{L/K}(\J_L)=\R_K N_{L/K}(\U_L^{\f_{\si{L}}})$, avec $\U_L^{\f_{\si{L}}}=\prod_{\p\nmid\si{\infty}}\U_{K_\p}^{v_\p(\f_{\si{K}})}\prod_{\p\mid\si{\infty},\,\p\nmid\f_{\si{K}}}\R_{K_\p}$.}\smallskip

\noindent En particulier, elle est non-ramifiée aux places finies étrangères à $\f_K$ et complètement composée aux places infinies étrangères à $\f_K$.
Si, maintenant, $\f_L$ est le diviseur de $L$ défini dans le corollaire, un calcul classique (cf. \cite{Ser}  montre que l'on a:\smallskip

\centerline{$N_{L/K}(\U_L^{\f_{\si{L}}})=\U_K^{\f_{\si{K}}}$; d'où: $\R_K N_{L/K}(\U_L^{\f_{\si{L}}})=\R_K\U_K^{\f_{\si{K}}}=\R_KN_{L/K}(\J_L)$.}\smallskip

\noindent Les conditions du Théorème \ref{Th12} sont donc remplies avec $\V_L=\R_L\U_L^{\f_{\si{L}}}$ dans l'extension abélienne $L/K$. Traduite en termes de classes de rayons, la conclusion $j_{K/L}(\J_K)\subset \R_L\U_L^{\f_{\si{L}}}$ fournit le résultat.

\subsection{Retour sur la conjecture principale de la capitulation}

Revenons maintenant sur la conjecture 4 et considérons pour cela le schéma d'extensions associé à une extension galoisienne non ramifiée $\si{\infty}$-décomposée $L/K$ de corps de nombres, où nous avons fait figurer les corps de classes de Hilbert respectifs de $L$ et de $K$:

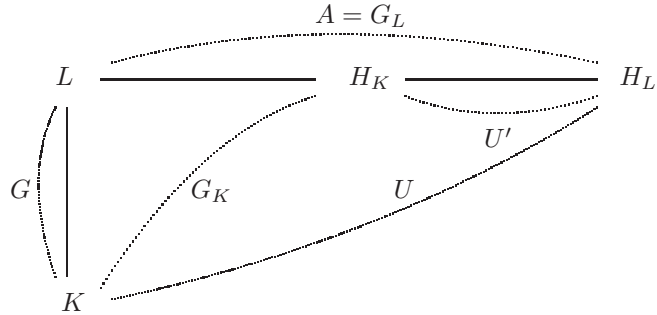
\begin{center}
\unitlength=1.5cm
\begin{picture}(6.6,2.8)

\put(0.65,0){$K$}
\put(0.7,0.3){\line(0,1){1.5}}
\put(0.6,2){$L$}

\bezier{60}(0.6,0.3)(0.3,1.2)(0.6,1.8)
\put(0.2,1.0){$G$}

\put(1.0,2.05){\line(1,0){1.9}}
\put(3.2,2){$H_K$}
\put(3.7,2.05){\line(1,0){1.7}}
\put(5.6,2){$H_L$}

\bezier{120}(1.1,2.2)(3,2.7)(5.4,2.2)
\put(2.9,2.55){$A=G_L$}

\bezier{70}(1.0,0.2)(1.9,1.6)(2.9,1.9)
\put(1.8,1.0){$G_K$}

\bezier{50}(3.7,1.9)(4.5,1.6)(5.4,1.9)
\put(4.4,1.45){$U'$}

\bezier{180}(1.1,0.1)(3.5,0.6)(5.4,1.8)
\put(3.6,1.0){$U$}
\end{picture}
\end{center}
\medskip\medskip

Le théorème de Tannaka-Terada ou plutôt le Corollaire \ref{C11} appliqué avec $U = \Gal(H_L/K)$ et $A = \Gal(H_K/L)$ nous montre que l'ordre de la capitulation $Cap_{L/K}$ dans l'extension $L/K$ est un multiple du degré $[L:K] = (U:A)$ pourvu qu'il existe un endomorphisme $\gamma$ du groupe $U$ tel que nous ayons:
$$
A / U' = (U/U')^{\gamma-1};
$$
de sorte que la Conjecture \ref{CP} est vérifiée dans ce cas.\smallskip

De cette condition K. Miyake donne la présentation cohomologique suivante (cf. \cite{Mi2}):
regardons le groupe $U$ comme extension du sous-groupe normal abélien $A$ par le quotient abélien $G = U/A$; faisons agir $U$ sur $A$ via les automorphismes intérieurs de $U$; et considérons l'homomorphisme $\pi$:\smallskip

\centerline{$H^1(U,A) \to H^1(U,A/U' )$,}\smallskip

\noindent induit par la surjection naturelle de $A$ sur $A/U'$ . Le groupe $U$ agissant trivialement sur le quotient $A/U'$ , nous avons évidemment:\smallskip

\centerline{$H^1(U , A / U' ) = \Hom (U,A/U') \simeq \Hom(U/U' ,A/U') = \Hom(G,A/U')$.}\smallskip

Pour tout élément $\varphi$ de $\Im \pi$, regardé comme homomorphisme de $G$ dans $A/U'$, posons:\smallskip

\centerline{$d_\varphi = |\Coker \varphi | = (A/U' : \Im\varphi)$.}\smallskip

\noindent Cela étant, nous avons:

\begin{Prop}\label{P15}
Étant donné un sous-groupe normal abélien $A$ d'un groupe $U$ qui définit un quotient abélien $U/A = G$, pour qu'il existe un endomorphisme surjectif $\gamma$ de $U$ tel qu'on ait:
$$
A = U^{\gamma-1}U',
$$
il est nécessaire et suffisant qu'il existe un morphisme surjectif $\varphi \in \Hom (U,A/U')$ qui soit la réduction modulo $U'$ d'un 1-cocycle $f$ de $U$ dans $A$. Lorsque cette condition est remplie, l'ordre $(U:A)$ du groupe quotient $U/A$ divise celui du noyau de l'homomorphisme de transfert $\Ver_{U/A}: U/U' \to A$.
\end{Prop}

\begin{Sco}
Plus généralement, pour tout homomorphisme $\varphi \in \Hom (U,A/U')$ qui provient d'un 1-cocycle $f \in Z^1(U,A)$, on a  l'identité:
$$
(\Ker\varphi)^{d_\varphi} \subset \Ker \Ver_{U/A},
$$
où $d_\varphi = (A/U' : \Im \varphi)$ est l' ordre du conoyau de $\varphi$.
\end{Sco}

\Preuve Comme exposé plus haut, tout 1-cocycle $f$ de $U$ dans $A$ induit, par réduction modulo $U'$, un morphisme $\varphi$ de $U/U'$ dans $A/U'$. Inversement à un tel cocycle $f$ correspond canoniquement un endomorphisme $\gamma$ de $U$, relevant $\varphi$, et défini par l'identité:\smallskip

\centerline{$\gamma ( x ) = f(x)x , \quad \forall x \in U$.}\smallskip

\noindent Le sous-groupe normal $A_\gamma$ qui lui est associé par la relation\smallskip

\centerline{$A_\gamma = U^{\gamma-1}U'$.}\smallskip

\noindent n'est autre que  l'image réciproque dans $A$ du sous groupe $\Im \varphi$ de $A/U'$. Il vient ainsi:\smallskip

\centerline{$d_\varphi = (A/U' : \Im \varphi) = (A:A_\gamma)$,}\smallskip

\noindent donc $\Ver_{U/A_\gamma} =\Ver_{A/A_\gamma} \circ \Ver_{U/A} = (\Ver_{U/A})^{d_\varphi}$, et l'identité annoncée
$ (\Ker \varphi)^{d_\varphi}\subset \Ker \Ver_{U/A}$ résulte alors de la Proposition \ref{P15} appliquée à $A_\gamma$.\smallskip

Peut-être n'est-il pas inutile de réécrire la Proposition \ref{P15} en terme de classes d'idéaux :

\begin{Cor}
Soient $L/K$ une extension abélienne non-ramifiée $\infty$-décomposée de corps de nombres algébriques, $H_K$ le corps des classes de Hilbert de $K$, et $H_L$ celui de $L$. La conjecture 4 est vérifiée dans $L/K$ sous la condition suffisante qu'il existe un morphisme surjectif $\varphi\in \Hom\big(\Gal (H_L/K) , N_{L/K}(Cl_L)\big)$ qui provienne d'un 1-cocycle $f \in Z^1\big(\Gal (H_L/K) ,Cl_L \big)$.
\end{Cor}

Tout comme dans la section précédente, le résultat obtenu est susceptible d'une interprétation idélique. Généralisant le Théorème \ref{Th12}, et avec les mêmes notations, nous obtenons, en effet:

\begin{Th}
Soient $K/k$ une extension galoisienne de corps de nombres et $G$ son groupe de Galois. Alors, si $\V_K$ est un sous-groupe ouvert de $\J_K$ contenant $\R_K$, qui est stable par $G$, on a:\smallskip

\centerline{$\big(j_{K/k}(\J_k) \V_K{}^{\si{-1}\!}N_{K/k}(\R_k)\big)^d \subset \big(j_{K/k}(\J_k)\cap \V_K\big)$, avec $d=\big(\J_K:(j_{K/k}(\J_k) \V_K{}^{\si{-1}\!}N_{K/k}(\R_k))\big)$.}
\end{Th}

\noindent {\em Preuve} (d'après Miyake \cite{Mi2}). Supposons d'abord que $K$ soit abélienne sur $k$ et introduisons l'extension abélienne $L_K$ de $K$ associée à $\V_K$ par la théorie du corps de classes. Notant $L_k$ la sous-extension maximale de $L_K$ qui est abélienne sur $k$, nous obtenons le schéma galoisien:

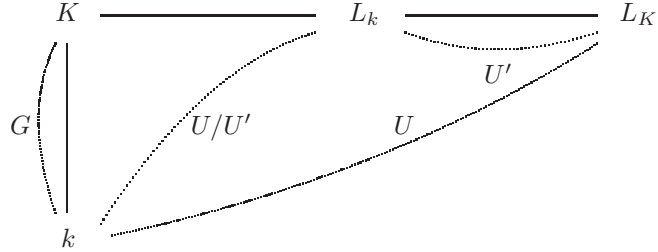
\begin{center}
\unitlength=1.5cm
\begin{picture}(6.6,2.8)

\put(0.65,0){$k$}
\put(0.7,0.3){\line(0,1){1.5}}
\put(0.6,2){$K$}

\bezier{60}(0.6,0.3)(0.3,1.2)(0.6,1.8)
\put(0.2,1.0){$G$}

\put(1.0,2.05){\line(1,0){1.9}}
\put(3.2,2){$L_k$}
\put(3.7,2.05){\line(1,0){1.7}}
\put(5.6,2){$L_K$}

\bezier{70}(1.0,0.2)(1.9,1.6)(2.9,1.9)
\put(1.8,1.0){$U/U'$}

\bezier{50}(3.7,1.9)(4.5,1.6)(5.4,1.9)
\put(4.4,1.45){$U'$}

\bezier{180}(1.1,0.1)(3.5,0.6)(5.4,1.8)
\put(3.6,1.0){$U$}
\end{picture}
\end{center}
\medskip\medskip

\noindent ainsi que les isomorphismes:\smallskip

\centerline{$A\simeq \J_K/V_K$;\qquad $U'\simeq \V_K{}^{\si{-1\!}}N_{K/k}(\R_k)/\V_K$; \qquad $U/U'\simeq\J_k/N_{K/k}(\V_K)\R_k$.}
\smallskip

Considérons l'application naturelle $\varphi\in\Hom (U/U' , A/U' )$ induite par l'extension des idèles de $k$ à $K$. Elle provient bien entendu de l'application induite par l'extension des idèles dans $\Hom (\J_k/\R_k N_{K/k}(\V_K),(\J_K/\V_K)^G)$. Mais, comme ce dernier groupe s'écrit encore $\Hom (U/U',A^G) \simeq \Hom(U,A^G) = H^1(U,A^G)$, et qu'il est trivialement contenu dans $H^1(U,A)$, nous pouvons appliquer la proposition 15 qui nous donne directement:\smallskip

\centerline{$(\Ker\varphi)^{d_\varphi}\subset\Ker\Ver_{U/A}$,}\smallskip

\noindent ce qui est précisément le résultat annoncé.\smallskip

Le cas galoisien s'en déduit sans difficulté, comme le Théorème \ref{Th12}, qui correspond au cas particulier $\V_K{}^{\si{-1}\!}N_{K/k}(\R_k)=\J_K$.

\section*{Bibliographie}
\addcontentsline{toc}{section}{Bibliographie}

Vu la complexité du sujet, nous avons choisi d'organiser exceptionnellement la bibliographie autour de trois thèmes directeurs, en respectant strictement dans chacune des sections la chronologie des publications.

\def\refname{\small{\sc Autour du Théorème 94 de Hilbert}}

\def\refname{\small{\sc Autour du Théorème d'Artin-Furtwängler}}

\def\refname{\small{\sc Autour du Théorème de Tannaka-Terada}}

\def\refname{\small{\sc Et aussi \dots}}

\section*{Addendum}
\addcontentsline{toc}{section}{Addendum}
Ce qui précède n'est autre que la mise au format \TeX, après correction de nombre d'erreurs typographiques, de la synthèse sur la capitulation des groupes de classes des corps de nombres publiée au Séminaire de Théorie des Nombres de Bordeaux et reçue le 14 octobre 1988.\medskip

Parmi les travaux ultérieurs, il convient de signaler le rapport de Miyake \cite{Mi6}, paru peu après, et, bien sûr, la démonstration algébrique de la conjecture principale donnée finalement par Suzuki \cite{Suz} deux ans plus tard. Une approche alternative débouchant effectivement sur une nouvelle preuve a été proposée plus récemment par Gruenberg et Weiss \cite{GW1}, lesquels l'ont généralisée dans \cite{GW2}.

\def\refname{\small{\sc Complément bibliographique}}

\bigskip
{\small
\noindent{\sc Adresse:}
Univ. Bordeaux \& CNRS,\\
Institut de Mathématiques de Bordeaux,\\
351 Cours de la Libération,
F-33405 Talence cedex

\noindent{\sc Courriel:}
 \tt jean-francois.jaulent@math.u-bordeaux.fr\\
 \url{https://www.math.u-bordeaux.fr/~jjaulent/}
}

\end{document}